\documentclass[10pt]{article}
\usepackage[utf8]{inputenc}
\usepackage{amsmath}
\usepackage{amsthm}
\usepackage{amsfonts}
\usepackage{amssymb}
\usepackage{color}
\usepackage[all]{xy}
\usepackage{pdfpages}
\usepackage{url}

\newtheorem{teo}{Theorem}[section]
\newtheorem{lema}{Lemma}[section]

\newtheorem{de}{Definition}[section]

\newtheorem{pr}{Proposition}[section]

\title{Classification of Toric 2-Fano 4-folds}
\author{Edilaine Ervilha Nobili}
\date{}

\begin{document}

\maketitle
 
\begin{abstract} In this notes we classify toric Fano 4-folds having positive second Chern Character. 
\end{abstract}
\textbf{keywords:} Chern character, toric variety, invariant surface, 2-Fano.

\section{Introduction}
 A smooth complex projective variety $X$ is said to be Fano if it has ample anti-canonical divisor. These varieties have been studied by several authors and play an important role in birational algebraic geometry. Fano varieties are quite rare. It was proved by Kollár, Miyaoka and Mori that, fixed the dimension, there exist only finitely many smooth Fano varieties up to deformation  (see \cite{kollar}, \cite{kol}). Further, in the toric case, there exist finitely many isomorphism classes of them. In dimension 2 there are five isomorphism classes of smooth toric Fano varieties  (also known as Del Pezzo surfaces): $\mathbb{P}^2$, $\mathbb{P}^1 \times \mathbb{P}^1$ and $\mathbb{P}^2$ blown up in 1, 2 or 3 points in general position. In dimensions 3, 4 and 5 there are 18, 124 and 866 isomorphism classes respectively (see \cite{bat3}, \cite{bat}, and \cite{nill}).\\
 
  In this paper we are interested in toric varieties known as 2-Fano varieties. A Fano variety $X$ is said to be 2-Fano if its second Chern character is positive (i.e., $ch_2(T_X)\cdot S >0$ for every surface $S\subset X$). These varieties were introduced by de Jong and Starr in \cite{starr} and \cite{jason} in connection with rationally simply connected varieties, which in turn are linked with the problem of finding rational sections for fibrations over surfaces. 2-Fano varieties are even rarer than Fano varieties. By (\cite{jason}, Theorem 1.2), the only 2-Fano surface is $\mathbb{P}^2$. In \cite{carol} it is proved that the only 2-Fano threefolds are $\mathbb{P}^3$ and the smooth hyperquadric in $\mathbb{P}^4$. In higher dimensions, few examples are known. First de Jong and Starr gave some examples in \cite{jason}, then Araujo and Castravet found some more examples (see \cite{ana}, section 5). Among all examples known, the only smooth toric 2-Fano varieties are projective spaces. 
  In this work, we go through the classification of toric Fano 4-folds, given by Batyrev in \cite{bat}, and we check that the only one with positive second Chern character is $\mathbb{P}^4$. We remark that in  \cite{sato2} Sato considers a similar problem. In particular, he classifies smooth toric Fano varieties with Picard number 2 whose second Chern character is nef (i.e., $ch_2(T_X)\cdot S \geq 0$ for every surface $S\subset X$).\\

\textbf{Notation:}	Throughout  the paper, we work over $\mathbb{C}$ and follow the standard notation for toric varieties.
Let $N$ be a lattice of rank $n$, $M:=N^{\vee}=Hom_\mathbb{Z}(N,\mathbb{Z})$, $N_\mathbb{R}:=N\otimes_\mathbb{Z}\mathbb{R}$ and $\Sigma\subseteq N_\mathbb{R}\simeq\mathbb{R}^n$ a simplicial complete fan. We will denote by $\Sigma(1)=\{v_1,...,v_d\}$ the set of all minimal generators of the fan $\Sigma$, $\Sigma(k)$ the set of k-dimensional cones of $\Sigma$ for $k>1$ and $X:=X_{\Sigma}$ the complete $\mathbb{Q}$-factorial toric variety defined by $\Sigma$. If $\sigma$ is a cone of $\Sigma$, the affine toric variety of the cone $\sigma$ will be denoted by $U_{\sigma}:=Spec(\mathbb{C}[S_{\sigma}])$, where $S_{\sigma}=\sigma^\vee \cap M$. We will denote by $V(\sigma)$ the closure of $Orbit(x_\sigma)$, where $x_\sigma$ is the distinguished point of $\sigma$. 

\section{Preliminaries}
\begin{de} A subset $\mathcal{P}=\{v_{i_1},v_{i_2},...,v_{i_r}\}$ of $\Sigma(1)$ is a \textbf{primitive collection} for $\Sigma$ if the following conditions are satisfied:\\
\begin{enumerate}
 \item 
$\mathcal{P}$ is not contained in any cone from $\Sigma$.
\item Any proper subset of $\mathcal{P}$ is contained in some cone from $\Sigma$.
 \end{enumerate}
 
\end{de}

\begin{de}  Let $\mathcal{P}=\{v_{i_1},v_{i_2},...,v_{i_r}\}\subseteq\Sigma(1)$ be a primitive collection for $\Sigma$ and $\sigma_{\mathcal{P}}=\langle{v_{j_1},...,v_{j_k}}\rangle$ be the cone of minimal dimension in $\Sigma$ such that $v_{i_1}+v_{i_2}+...+v_{i_r}\in\sigma_{\mathcal{P}}$. Then, there is a unique linear relation $$v_{i_1}+...+v_{i_r}=c_1v_{j_1}+...+c_kv_{j_k}, \ \ \ c_i\in\mathbb{Q}_{>0}.$$
We call  $\mathcal{R(P)}:=v_{i_1}+...+v_{i_r}-c_1v_{j_1}-...-c_kv_{j_k}$ the \textbf{primitive relation} associated to $\mathcal{P}$.

If $X$ is smooth then $c_i\in\mathbb{Z}_{>0}$, and we define the \textbf{degree} of the primitive collection $\mathcal{P}$ by:
$$\Delta(\mathcal{P}):=r-c_1-...-c_k.$$
\end{de}

Let $A_1(X)$ be the group of 1-cycles on $X$ modulo numerical equivalence, $N_1(X)=A_1(X)\otimes_{\mathbb{Z}}\mathbb{R}$, and $NE(X)\subseteq N_1(X)$ the Mori cone of $X$, generated by classes of irreducible complete curves. It is well known that $NE(X)=\displaystyle\sum_{\sigma \in \Sigma(n-1)}\mathbb{R}_{\geq0}[V(\sigma)]$ (see for instance \cite{cox}, Theorem 6.2.20). 
We also recall that there is an exact sequence:

$$
\begin{array}{ccccccccc}
 0 &\longrightarrow& N_1(X) &\longrightarrow& \mathbb{R}^d&\longrightarrow&  N_{\mathbb{R}} &\longrightarrow &0 \ \ ,\\
  &&\xi&\longmapsto &\displaystyle\sum_{v_i\in\Sigma(1)}(V(\langle v_i \rangle)\cdot \xi)e_i  \\
&&&&e_i&\longmapsto & v_i \\
\end{array}
$$

 \noindent where $d:=\sharp \ \Sigma(1)$, and ${e_1,...,e_d}$ is the standard basis for $\mathbb{R}^d$.\\
\\
In particular, we may interpret $N_1(X)$ as the space of linear relations among the minimal generators of $\Sigma$ and, under this identification, we have that $N_1(X)$ is generated by primitive relations. Moreover, if $V(\sigma)$ is an invariant curve on $X$ then its class $[V(\sigma)]$ in $N_1(X)$ is a positive linear combination of primitive relations (see \cite{cox}, Theorem 6.3.10).\\

Hence, \  $NE(X)=\displaystyle\sum_{\mathcal{P} \ \substack{ primitive \\  collection}}\mathbb{R}_{\geq0}\mathcal{R(P)}$.

Note that a relation $\displaystyle\sum_{i=1}^{d}a_iv_i=0, a_i \in \mathbb{R}$, corresponds to an element $\xi\in N_1(X)$ that has intersection $a_i$ with $V(\langle v_i \rangle)$ for all $i\in \{1,...,d\}$.\\
Since we known $c_1(T_X)=\displaystyle \sum_{i=1}^{d}V(\langle v_i \rangle)$ (see for instance, \cite{cox} Theorem 8.2.3), if $X$ is smooth and $\mathcal{R(P)}=\displaystyle\sum_{i=1}^{d}a_iv_i$ is a primitive relation, then $$\Delta(\mathcal{P})=\displaystyle\sum_{i=1}^{d}a_i=\displaystyle\sum_{i=1}^{d}V(\langle v_i \rangle)\cdot\mathcal{R(P)}=-K_X\cdot\mathcal{R(P)}.$$
Hence, using Kleiman's Criterion of ampleness, we can give a characterization of smooth toric Fano varieties in terms of primitive relations:\\

A smooth toric variety $X_{\Sigma}$ is a Fano variety  if and only if $\Delta(\mathcal{P})>0$ for every primitive relation $\mathcal{P}$ of $\Sigma$.\\

From now on, $X:=X_{\Sigma}$ will denote a smooth projective toric variety.\\
In (\cite{bat}, 2.2.4) we see that Fano toric varieties can be recovered from the set of primitive relations. In that paper Batyrev gives a classification of toric Fano 4-folds by describing the possible sets of primitive relations in dimension 4. He also gives
  a geometric description for these varieties. He found 123 isomorphism classes of smooth toric Fano 4-fold. Then in \cite{sato} Sato noticed one missing isomorphism class in Batyrev's classification and he described the primitive relations of this missing class, completing the classification of toric Fano 4-folds. \\
 Note that if $X_{\Sigma}$ is a smooth toric variety then any set of primitive vectors that generate a maximal cone of $\Sigma$ can be chosen to be the canonical basis of $\mathbb{Z}^n$. By definition of primitive collection, a cone $\sigma$ belongs to $\Sigma$ if and only if $\sigma$ does not contain any primitive collection. In the next example, we illustrate how to recover a smooth toric Fano variety from its set of primitive relations.
\\

  \textbf{Example}($H_1$):
Let $X_{\Sigma}$ be the toric Fano 4-fold given by the following primitive relations:\\
$v_1+v_2=v_8, \ v_7+v_8=v_1, \ v_1+v_6=v_7,  \ v_2+v_7=0, \ v_6+v_8=0, \ v_3+v_4+v_5=2v_1$. 

Then, the primitive collections are: $\{v_1,v_2\}$, $\{v_7,v_8\}$, $\{v_1,v_6\}$, $\{v_2,v_7\}$, $\{v_6,v_8\}$, $\{v_3,v_4,v_5\}$.
Thus, by definition of primitive collection, the fan $\Sigma$ obtained from primitive relations, satisfies:\\
$$\sigma=\langle v_i,v_j,v_k,v_l\rangle \in \Sigma \Leftrightarrow
\left\{
\begin{array}{lll}
\langle v_1,v_2\rangle & \nsubseteq & \sigma \\
\langle v_7,v_8\rangle & \nsubseteq & \sigma \\
\langle v_1,v_6\rangle & \nsubseteq & \sigma \\
\langle v_2,v_7\rangle & \nsubseteq & \sigma \\
\langle v_6,v_8\rangle & \nsubseteq & \sigma \\
\langle v_3,v_4,v_5\rangle & \nsubseteq & \sigma \\
\end{array}
\right.
$$
Since $X_{\Sigma}$ is smooth, every maximal cone in $\Sigma$ provides a basis to $N\simeq \mathbb{Z}^4$.
The cone $\langle v_1,v_2,v_3,v_4\rangle$ is maximal in $\Sigma$, so we can take $v_1=(1,0,0,0), v_2=(0,1,0,0), v_3=(0,0,1,0), v_4=(0,0,0,1)$. Thus, from the primitive relations, we get $v_5=(2,0,-1,-1), v_6=(-1,-1,0,0), v_7=(0,-1,0,0), v_8=(1,1,0,0)$.
  \\
 
 In the table below we list all smooth toric 4-folds and its primitive collections or geometric description. The last variety in our table follows Sato's notation. For the others, our notation differs from Batyrev's notation used in \cite{bat} only in the enumeration of  minimal vectors. Whenever he enumerates the vectors from $0$ to $k$, we will enumerate them from $1$ to $k+1$. We denote by $S_i$ the Del Pezzo surface obtained by the blow up of $i$ points in general position on $\mathbb{P}^2$ for $i=1,2$ and 3. It is clear that primitive collections are not enough to describe the variety. They describe only its combinatorial type.  \\ \\

 \begin{tabular}{|p{2cm}| p{13,5cm}|}
\hline
Notation & \hspace{2.5cm} Primitive Collections or Geometric Description  \\ 

\hline
& $\mathbb{P}^4$\\ \hline

$B_1,...,B_5$ & $\{v_5,v_6\}$, $\{v_1, v_2, v_3, v_4\}$ \\ \hline

$C_1,...,C_4$ & $\{v_1, v_2, v_3\}$, $\{v_4, v_5, v_6\}$\\ \hline

$D_1,...,D_{19}$ &  $\{v_4,v_5\}$, $\{v_6, v_7\}$, $\{v_1, v_2, v_3\}$ \\ 
\hline

$E_1,E_2,E_3$ & $\{v_1,v_7\}$, $\{v_1, v_2\}$, $\{v_6, v_7\}$, $\{v_2, v_3, v_4, v_5\}$, $\{v_3, v_4, v_5, v_6\}$  \\
\hline

 $G_1,...,G_6$ & $\{v_1, v_7\}$, $\{v_2, v_3, v_4\}$, $\{v_4, v_5, v_6\}$, $\{v_5, v_6, v_7\}$, $\{v_1, v_2, v_3\}$ \\
\hline

$H_1,...,H_{10}$ & $\{v_1, v_2\}$, $\{v_7, v_8\}$, $\{v_1, v_6\}$, $\{v_2, v_7\}$, $\{v_6, v_8\}$, $\{v_3, v_4, v_5\}$\\
\hline

 $I_1,...,I_{15}$ & $\{v_1,v_2\}$, $\{v_7, v_8\}$, $\{v_3,v_6\}$, $\{v_6,v_8\}$, $\{v_3, v_4, v_5\}$, $\{v_4, v_5, v_7\}$ \\
\hline
\end{tabular}

\begin{tabular}{|p{2cm}| p{13,5cm}|}
\hline
$J_1,J_2$ & $\{v_3,v_6\}$, $\{v_6, v_8\}$, $\{v_7,v_8\}$, $\{v_1, v_2, v_3\}$, $\{v_1, v_2,v_7\}$, $\{v_1,v_2,v_8\}$, $\{v_3, v_4, v_5\}$, $\{v_4, v_5, v_6\}$, $\{v_4, v_5, v_7\}$ \\
\hline

$K_1,...,K_4$ & $\{v_7, v_9\}$, $\{v_1, v_8\}$, $\{v_8, v_9\}$, $\{v_2, v_8\}$, $\{v_6, v_7\}$, $\{v_1, v_6\}$, $\{v_6, v_9\}$, $\{v_1, v_2\}$, $\{ v_2, v_7\}$, $\{v_3, v_4, v_5\}$ \\
\hline

$L_1,...,L_{13}$ & $ \{v_1, v_8\}$, $\{v_2, v_3\}$, $\{v_4, v_5\}$, $\{v_6, v_7\}$  \\ \hline

$M_1,...,M_4$ & $\{v_1, v_8\}$, $\{v_4, v_5\}$, $\{v_6, v_7\}$, $\{v_1, v_2, v_3\}$, $\{v_4, v_6, v_8\}$, $\{v_2, v_3, v_5\}$, $\{v_2, v_3, v_7\}$  \\
\hline

$Q_1,...,Q_{17}$ & $\{v_1,v_2\}$, $\{v_1, v_8\}$, $\{v_2,v_7\}$, $\{v_3, v_5\}$, $\{v_4, v_6\}$, $\{v_8,v_9\}$, $\{v_7, v_9\}$ \\
\hline

$R_1,R_2,R_3$ & $\{v_7,v_9\}$, $\{v_4,v_8\}$, $\{v_8,v_9\}$, $\{v_6,v_7\}$, $\{v_3,v_5\}$, $\{v_4,v_6\}$ $\{v_1, v_2, v_9\}$, $\{v_3, v_6, v_8\}$, $\{v_1, v_2, v_5\}$, $\{v_1, v_2, v_7\}$, $\{v_1, v_2, v_4\}$ \\
\hline

$108$ & $\{v_7,v_9\}$, $\{v_8, v_9\}$, $\{v_3,v_5\}$, $\{v_4, v_6\}$, $\{v_1, v_7\}$, $\{v_3,v_6\}$, $\{v_1, v_2, v_5\}$, $\{v_1, v_2, v_4\}$, $\{v_2, v_5, v_8\}$, $\{v_2, v_4, v_8\}$ \\ \hline 

$U_1,...,U_8$ & $\{v_1, v_3\}$, $\{v_2, v_4\}$, $\{v_1, v_4\}$, $\{v_3, v_5\}$, $\{v_4, v_6\}$, $\{v_2, v_5\}$, $\{v_1, v_5\}$, $\{v_2, v_6\}$, $\{ v_3, v_6\}$,  $\{ v_7, v_8\}$,  $\{ v_9, v_{10}\}$ \\
\hline

$Z_1,Z_2$ & $\{v_1, v_8\}$, $\{v_5, v_7\}$, $\{v_1, v_2, v_5\}$, $\{v_1, v_2, v_6\}$, $\{v_2, v_4, v_5\}$, $\{v_2, v_4, v_6\}$, $\{v_3, v_7, v_8\}$, $\{v_3, v_4, v_6\}$, $\{v_3, v_4, v_7\}$, $\{v_3, v_6, v_8\}$   \\
\hline
$117$ & $\{v_4, v_{10}\}$, $\{v_1, v_5\}$, $\{v_2, v_6\}$, $\{v_3, v_7\}$, $\{v_8, v_9\},$ $\{v_1, v_2, v_{10}\}$, $\{v_1, v_3, v_{10}\}$, $\{v_2, v_3, v_{10}\}$, $\{v_1, v_2, v_3\}$, $\{v_1, v_9, v_{10}\}$, $\{v_2, v_9, v_{10}\}$, $\{v_3, v_9, v_{10}\},$  $\{v_1, v_2, v_9\}$, $\{v_1, v_3, v_9\}$, $\{v_2, v_3, v_9\}$, $\{v_4, v_5, v_6\}$, $\{v_4, v_5, v_7\}$, $\{v_4, v_6, v_7\}$, $\{v_5, v_6, v_7\}$, $\{v_4, v_5, v_8\}$, $\{v_4, v_6, v_8\}$, $\{v_4, v_7, v_8\}$, $\{v_5, v_6, v_8\}$, $\{v_5, v_7, v_8\}$,  $\{v_6, v_7, v_8\}$           \\
\hline
$118$ & $\{v_4, v_9\}$,  $\{v_1, v_5\}$, $\{v_2, v_6\}$, $\{v_3, v_7\}$,  $\{v_1,v_2,v_9\}$, $\{v_1, v_3, v_9\}$, $\{v_2, v_3, v_9\}$,  $\{v_1, v_2, v_3\}$, $\{v_4, v_5, v_8\}$, $\{v_4, v_6, v_8\}$, $\{v_4, v_7, v_8\}$, $\{v_5, v_6, v_8\}$, $\{v_5, v_7, v_8\}$, $\{v_6, v_7, v_8\}$   \\ \hline 

$119, 120, 121$ & $S_2\times S_2, \ S_2\times S_3, \ S_3\times S_3$ \\ \hline
$124$ & $\{v_1, v_4\}$,  $\{v_2, v_5\}$, $\{v_3, v_6\}$, $\{v_1, v_2, v_3\}$,  $\{v_4,v_5,v_6\}$, $\{v_7, v_8, v_9\}$,  $\{v_1, v_2, v_9\}$, $\{v_4, v_5, v_9\}$, $\{v_1, v_3, v_8\}$, $\{v_4, v_6, v_8\}$, $\{v_2, v_3, v_7\}$, $\{v_5, v_6, v_7\}$, $\{v_1, v_8, v_9\},$ $\{v_4, v_8, v_9\}$, $\{v_2, v_7, v_9\}$, $\{v_5, v_7, v_9\}$, $\{v_3, v_7, v_8\}$, $\{v_6, v_7, v_8\}$ \\ \hline

\end{tabular}\\ \\

\textbf{Remark}. There is a misprint in (\cite{bat}, Proposition 3.4.1) concerning the primitive relations for the toric Fano 4-fold 108.

\section{Second Chern Class Computation}

We want to use Batyrev's classification of toric Fano 4-folds to classify toric 2-Fano 4-folds, that is, smooth toric Fano 4-folds with positive second Chern character. For this purpose, we will compute $ch_2(T_X)$ in terms of the invariant divisors $D_i:=V(\langle v_i\rangle)$.\\

\begin{pr} For a smooth toric variety $X$ we have: $$ ch_2(T_X)=\frac{1}{2}\left(\displaystyle \sum_{i=1}^{d}D_i^2\right).$$ 
\end{pr}

\textit{Proof.} There are exact sequences: (\cite{cox}, 4.0.28, 8.1.1)\\
$$0\rightarrow\Omega_X^1\rightarrow \mathcal{O}_X^n\rightarrow \displaystyle \oplus_{i=1}^{d}\mathcal{O}_{D_i}\rightarrow0$$
$$0\rightarrow\mathcal{O}(-D_i)\rightarrow \mathcal{O}_X\rightarrow \mathcal{O}_{D_i}\rightarrow0$$ 
Where $\mathcal{O}_{D_i}$ is the structure sheaf on $D_i$ extended by zero to X.\\

Using Whitney sum we have:\\

$0=ch_2(\mathcal{O}_{X{_\Sigma}}^n)=ch_2(\Omega_X^1)+ch_2(\displaystyle \oplus_{i=1}^{d}\mathcal{O}_{D_i})$\\

$ch_2(\mathcal{O}_{D_i})=-ch_2(\mathcal{O}(-D_i)=-\frac{1}{2}D_i^2$ \ for all $i=1,...,d$.\\

$\Rightarrow ch_2(T_X)=ch_2(\Omega_X^1)=-ch_2(\displaystyle \oplus_{i=1}^{d}\mathcal{O}_{D_i})=\frac{1}{2}\left(\displaystyle \sum_{i=1}^{d}D_i^2\right)$. \qed\\

By definition, a variety $X$ has positive second Chern character if for any surface $S\in X$ we have $ch_2(X)\cdot S>0$. However, in the toric case, we only need to check this inequalities for invariant surfaces, because of the following result. The proof sketched below is due to D. Monsôres.

\begin{pr}Let $X:= X_{\Sigma}$ be a complete toric variety of dimension $n \geq 3$. If $S$ is a surface
contained in $X$, then we have a numerical equivalence: 

\begin{center}     $S \equiv \displaystyle\sum_{\sigma \in \Sigma(n-2)} a_{\sigma} \cdot
[V(\sigma)]$

\end{center}

 \noindent with $a_{\sigma} \geq 0$, \  $\forall \sigma \in \Sigma(n-2)$.

\end{pr}

\textit{Sketch of the proof}. The proof is by induction on the dimension of $X$. If $n=3$ then $S$ is an efective divisor on $X$ and by (\cite{fulton}, section 5.1) $S\sim \displaystyle\sum_{\sigma\in \Sigma(1)}a_{\sigma}\cdot [V(\sigma)]$ with $a_{\sigma}\geq 0, \forall \sigma\in \Sigma(1).$ By induction hypothesis we can suppose that $S$ intersects the torus $T$ of $X$. Consider the action $\mathbb{C}^*\times X\rightarrow X$ given by $(t,x)\mapsto t^{\lambda_1}\cdot x$, where $\lambda_1\in N$. This action induces a rational map $f:\mathbb{C}\times S\dashrightarrow X$. Consider a toric resolution on indeterminacy for this map:\\

$$\xymatrix{
Y \ar@/_0.8cm/[dd]_{\pi}
 \ar[d]_p \ar[dr]^{\psi} & \\ 
\mathbb{C}\times S \ar[d]_q \ar@{-->}[r]^f  & X  \\
\mathbb{C} & & 
}$$
 By (\cite{hartshorne},III 9.6,9.7) $\pi$ is a flat morphism whose fibers have pure dimension two. Hence the cycles $\psi_*(\pi^*(0))$ and $\psi_*(\pi^*(1))$ are rationally equivalent. Since $S=\psi_*(\pi^*(1))$ we get a numerical equivalence $S\equiv \displaystyle\sum_{i=1}^ {k}a_iS_i$, where $a_i\geq 0$ and $S_i\subset \psi(\pi^{-1}(0)) \ \forall i=1,...,k$. Note that by construction $\psi_*(\pi^*(0))$ and therefore every $S_i$ is invariant by the action of $\lambda_1$. If each surface $S_i$ has empty intersection with the torus then we conclude the proposition by induction. If $S_i$ intersects the torus we take $\lambda_2\in N$ such that $\{\lambda_1, \lambda_2\}$ are linearly independent and repeat the construction above to $S_i$ and $\lambda_2$. We get $S_i\equiv \displaystyle\sum_{j=1}^{r}b_jS'_j$, where $b_j\geq 0$ and each $S'_j$ is an invariant surface by the actions of $\lambda_1$ and $\lambda_2$. If $S'_j\cap T=\emptyset \ \forall j=1,...,r$ we are done. If for some $S'_j$ it fails, we repeat the process using a parameter $\lambda_3$ such that $\{\lambda_1,\lambda_2,\lambda_3\}$ are linearly independent and obtain $S'_j\equiv \displaystyle\sum_{k=1}^{s}c_kS''_k$ where $c_k\geq 0$ and $S''_k$ are invariant surfaces by the actions of $\lambda_1, \lambda_2$ and $\lambda_3$. To finish the proof we observe that $S''_k\cap T=\emptyset \ \forall k=1,...,s$. Since $\{\lambda_1, \lambda_2,\lambda_3\}$ are linearly independent, if there was $t\in S''_k\cap T$ we would have an injective map $(\mathbb{C}^*)^3\rightarrow S''_k$ given by $(t_1,t_2,t_3)\mapsto t_1^{\lambda_1}\cdot t_2^{\lambda_2}\cdot t_3^{\lambda_3}\cdot t$. But this is absurd since $S''_k$ is a surface. \qed   \\

So, in order to check whether a smooth toric Fano variety is 2-Fano, we need to compute $\displaystyle \sum_{i=1}^{d}D_i^2 \cdot S$ for invariant surfaces $V(\sigma)$. For that purpose we will use the following formula (see for instance \cite{fulton}, 5.1).\\

\begin{lema} Let  $X_{\Sigma}$ be a smooth projective toric variety, $\sigma$ a cone in $\Sigma$ and $D_i=V(\langle v_i \rangle)$ an invariant divisor on $X_{\Sigma}$. Then, if $V(\sigma)\nsubseteq supp (D_i)$, we have:\\
 $$D_i \cdot V(\sigma)  =
\left\{
\begin{array}{ll}
  V(\langle v_i,\sigma \rangle)  & \text{if} \ v_i \ \text{and} \ \sigma \ \text{span a cone of}\  \Sigma\\
  0 & \text{if} \ v_i \ and \ \sigma \ \text{do not span a cone of}\  \Sigma \ \ . \\
\end{array}
\right.
$$
\\
\end{lema}
 
Suppose, otherwise, that $V(\sigma)\subseteq supp(D_i)$. Since $D_i$ is a Cartier divisor, there exists $u\in M$ such that $(D_i)_{\mid_{U_{\sigma}}}=div (\chi^ u)_{\mid_{U_{\sigma}}}$. Hence,  $D_i-div(\chi^u)$ is linearly equivalent to  $D_i$ and the support of $D_i - div(\chi^u)=D_i - \displaystyle \sum_{i=1}^{r}\langle u,v_j\rangle D_j$ does not contain $V(\sigma)$. Thus, if we find an element $u\in M$ satisfying $(D_i)_{\mid_{U_{\sigma}}}=div (\chi^ u)_{\mid_{U_{\sigma}}}$ then we can use the previous lemma to compute $D_i \cdot V(\sigma)=(D_i-div(\chi^ u)) \cdot V(\sigma)$.\\
By the cone-orbit correspondence, we have:
$$U_{\sigma}\cap D_j \neq \emptyset \Leftrightarrow \langle v_j\rangle \subseteq \sigma.$$
Since $div (\chi^ u)_{\mid_{U_{\sigma}}}=\displaystyle \sum_{v_j \in \sigma}\langle u,v_j\rangle (D_j)_{\mid_{U_{\sigma}}}$, in order to have $(D_i)\mid_{U_{\sigma}}=div(\chi^u)\mid_{U_{\sigma}}$ we can take $u$ to be any element in $M$ such that $\langle u,v_i\rangle =1$ and $\langle u,v_j\rangle =0 \ \forall j \neq i$ such that $v_j \in \sigma$. With this, we are ready to compute the product of $ch_2(T_X)$ with $V(\sigma)$.\\

\section{The Main Result}

\begin{teo} The only toric Fano 4-fold with positive second Chern character is $\mathbb{P}^4$.
\end{teo}

\textit{Proof}. First of all, it was shown in (\cite{jason}, Theorem 1.2) that in the following cases $ch_2(T_X)$ is not positive: 

\begin{enumerate}
\item $X=Z\times Y$ is a product of positive-dimensional Fano manifolds.
\item $X=\mathbb{P}_Y(E)$ is a projective bundle  over a positive-dimensional Fano manifold $Y$.
\end{enumerate}

As consequence, the toric Fano 4-folds listed in the Batyrev's classification that are of type 1. or 2. do not have positive second Chern Character. They are (see Batyrev's description of these varieties in \cite{bat}):\\
 
 $B_1,...,B_5$,$C_1,...,C_4,D_1,...,D_{19},H_8, L_1,...,L_{13},I_7,I_{11},I_{13},Q_6,Q_8,Q_{10}$,
 $Q_{11}$, $Q_{15}$,$K_4,U_4,U_5,U_6,119,120,121$. \\
 
In the remaining cases, we computed $ch_2(T_X)\cdot S$ for all invariant surfaces $S\subset X$, as described in section 3. To make the computation we used the program Maple. For all smooth toric Fano 4-folds $X\neq \mathbb{P}^4$ in Batyrev's list we found a surface $S\subset X$ such that $ch_2(T_X)\cdot S \leq 0$.\\

The next table summarizes our results. The first column lists toric Fano 4-folds according Batyrev's notation. The second column lists its primitive vectors explicity. The third column gives an invariant surface $S$ for which the intersection number $ch_2(T_X)\cdot S$ (listed on the last column) is non positive.\\ \\

\begin{tabular}{  |c|  p{11,2cm} | c| c| }
\hline
 & \centering{Primitive \ Vectors} & Surface & $ch_2(T_X)\cdot S$  \\ 
\hline
$E_1$ & $v_1=e_1, v_2=e_2,v_3=e_3,v_4=e_4,v_5=2e_1-e_2-e_3-e_4,v_6=e_1+e_2,v_7=-e_1$ & $V(v_2,v_3)$ &  -2 \\ \hline

$E_2$ & $v_1=e_1, v_2=e_2,v_3=e_3,v_4=e_4,v_5=e_1-e_2-e_3-e_4,v_6=e_1+e_2,v_7=-e_1$ & $V(v_2,v_3)$ & $-\displaystyle\frac{3}{2}$ \\ \hline

$E_3$ & $v_1=e_1, v_2=e_2,v_3=e_3,v_4=e_4,v_5=-e_2-e_3-e_4,v_6=e_1+e_2,v_7=-e_1$ & $V(v_2,v_3)$ & -1 \\ \hline

$G_1$ & $v_1=e_1, v_2=e_2,v_3=e_3,v_4=e_1-e_2-e_3,v_5=e_4,v_6=e_1+e_2+e_3-e_4,v_7=-e_1$ & $V(v_1,v_5)$ & $-\displaystyle\frac{1}{2}$ \\ \hline

$G_2$ & $v_1=e_1, v_2=e_2,v_3=-e_1-e_2,v_4=e_4,v_5=e_3,v_6=2e_1-e_3-e_4,v_7=-e_1+e_4$ & $V(v_1,v_5)$ & -2 \\ \hline

$G_3$ & $v_1=e_1, v_2=e_2,v_3=e_3,v_4=-e_2-e_3,v_5=e_4,v_6=e_1+e_2+e_3-e_4,v_7=-e_1$ & $V(v_1,v_5)$ & -1 \\ \hline

$G_4$ & $v_1=e_1, v_2=e_2,v_3=-e_1-e_2,v_4=e_4,v_5=e_3,v_6=e_1+e_2-e_3-e_4,v_7=-e_1+e_4$ & $V(v_1,v_5)$ & $-\displaystyle\frac{1}{2}$ \\ \hline

$G_5$ & $v_1=e_1, v_2=e_2,v_3=-e_1-e_2,v_4=e_4,v_5=e_3,v_6=-e_3-e_4,v_7=-e_1+e_4$ & $V(v_2,v_5)$ & -2 \\ \hline

$G_6$ & $v_1=e_1, v_2=e_2,v_3=-e_1-e_2,v_4=e_4,v_5=e_3,v_6=e_1-e_3-e_4,v_7=-e_1+e_4$ & $V(v_2,v_5)$ & $-\displaystyle\frac{3}{2}$ \\ \hline

$H_1$ & $v_1=e_1, v_2=e_2,v_3=e_3,v_4=e_4,v_5=2e_1-e_3-e_4,v_6=-e_1-e_2,v_7=-e_2,v_8=e_1+e_2$ & $V(v_3,v_4)$ & $-\displaystyle\frac{3}{2}$ \\ \hline
\end{tabular}

\begin{tabular}{  |c| p{11,2cm} | c| c| }
\hline
& \centering{Primitive \ Vectors} & Surface & $ch_2(T_X)\cdot S$  \\ 
\hline
$H_2$ & $v_1=e_1, v_2=e_2,v_3=e_3,v_4=e_4,v_5=2e_1+e_2-e_3-e_4,v_6=-e_1-e_2,v_7=-e_2,v_8=e_1+e_2$ & $V(v_3,v_4)$ & -1 \\ \hline

$H_3$ & $v_1=e_1, v_2=e_2,v_3=e_3,v_4=e_4,v_5=2e_1+2e_2-e_3-e_4,v_6=-e_1-e_2,v_7=-e_2,v_8=e_1+e_2$ & $V(v_3,v_4)$ & $-\displaystyle\frac{3}{2}$ \\ \hline

$H_4$ & $v_1=e_1, v_2=e_2,v_3=e_3,v_4=e_4,v_5=e_1-e_3-e_4,v_6=-e_1-e_2,v_7=-e_2,v_8=e_1+e_2$ & $V(v_3,v_4)$ & $-\displaystyle\frac{3}{2}$ \\ \hline

$H_5$ & $v_1=e_1, v_2=e_2,v_3=e_3,v_4=e_4,v_5=e_1+e_2-e_3-e_4,v_6=-e_1-e_2,v_7=-e_2,v_8=e_1+e_2$ & $V(v_3,v_4)$ & $-\displaystyle\frac{3}{2}$\\ \hline

$H_6$ & $v_1=e_1, v_2=e_2,v_3=e_3,v_4=e_4,v_5=e_1+2e_2-e_3-e_4,v_6=-e_1-e_2,v_7=-e_2,v_8=e_1+e_2$ & $V(v_3,v_4)$ & $-\displaystyle\frac{3}{2}$ \\ \hline

$H_7$ & $v_1=e_1$, $v_2=e_2$, $v_3=e_3$, $v_4=e_4$, $v_5=2e_2-e_3-e_4,$ $v_6=-e_1-e_2,$ $v_7=-e_2,$ $v_8=e_1+e_2$ & $V(v_3,v_4)$ & $-\displaystyle\frac{3}{2}$ \\ \hline

$H_9$ & $v_1=e_1, v_2=e_2,v_3=e_3,v_4=e_4,v_5=e_2-e_3-e_4,v_6=-e_1-e_2,v_7=-e_2,v_8=e_1+e_2$ & $V(v_3,v_4)$ & $-\displaystyle\frac{3}{2}$\\ \hline

$H_{10}$ & $v_1=e_1, v_2=e_2,v_3=e_3,v_4=e_4,v_5=-e_1-e_3-e_4,v_6=-e_1-e_2,v_7=-e_2,v_8=e_1+e_2$ & $V(v_3,v_4)$ & $-\displaystyle\frac{3}{2}$ \\ \hline

$I_1$ & $v_1=e_1, v_2=-e_1+e_3,v_3=e_3,v_4=e_4,v_5=-2e_2+e_3-e_4,v_6=e_2-e_3,v_7=e_2,v_8=-e_2+e_3$ & $V(v_1,v_4)$ & $-\displaystyle\frac{3}{2}$ \\ \hline

$I_2$ & $v_1=e_1, v_2=e_2,v_3=e_3,v_4=e_4,v_5=2e_1+2e_2-e_3-e_4,v_6=-e_1-e_2,v_7=-e_1-e_2+e_3,v_8=e_1+e_2$ & $V(v_1,v_4)$ & $-\displaystyle\frac{3}{2}$\\ 
\hline
$I_3$ & $v_1=e_1, v_2=e_2,v_3=e_3,v_4=e_4,v_5=2e_1+e_2-e_3-e_4,v_6=-e_1-e_2,v_7=-e_1-e_2+e_3,v_8=e_1+e_2$ & $V(v_1,v_4)$ & $-\displaystyle\frac{3}{2}$ \\ \hline

$I_4$ & $v_1=e_1, v_2=e_2,v_3=e_3,v_4=e_4,v_5=-2e_1-2e_2+e_3-e_4,v_6=e_1+e_2-e_3,v_7=e_1+e_2,v_8=-e_1-e_2+e_3$ & $V(v_1,v_4)$ & $-\displaystyle\frac{3}{2}$ \\ \hline
$I_5$ & $v_1=e_1, v_2=e_2,v_3=e_3,v_4=e_1+e_2,v_5=-e_1-e_2-e_3+2e_4,v_6=-e_4,v_7=e_3-e_4,v_8=e_4$ & $V(v_1,v_4)$ & $-\displaystyle\frac{3}{2}$ \\ \hline

$I_6$ & $v_1=e_1, v_2=e_2,v_3=e_1+e_2,v_4=e_4,v_5=e_3,v_6=-e_1-e_2-e_3-e_4,v_7=-e_3-e_4,v_8=e_1+e_2+e_3+e_4$ & $V(v_1,v_4)$ & $-\displaystyle\frac{3}{2}$ \\ \hline

$I_8$ & $v_1=e_1, v_2=e_2,v_3=e_3,v_4=e_4,v_5=-e_2-e_3-e_4,v_6=e_1+e_2,v_7=e_1+e_2+e_3,v_8=-e_1-e_2$ & $V(v_1,v_4)$ & $-\displaystyle\frac{3}{2}$ \\ \hline

$I_9$ & $v_1=e_1, v_2=e_2,v_3=e_3,v_4=e_4,v_5=-e_1-e_2-e_4,v_6=e_1+e_2-e_3,v_7=e_1+e_2,v_8=-e_1-e_2+e_3$ & $V(v_1,v_4)$ & $-\displaystyle\frac{3}{2}$\\ \hline

$I_{10}$ & $v_1=e_1, v_2=e_2,v_3=e_3,v_4=e_4,v_5=e_1+e_2-e_3-e_4,v_6=-e_1-e_2,v_7=-e_1-e_2+e_3,v_8=e_1+e_2$ & $V(v_1,v_4)$ & $-\displaystyle\frac{3}{2}$ \\ \hline

$I_{12}$ & $v_1=e_1, v_2=e_2,v_3=e_3,v_4=e_4,v_5=-e_1-e_2-e_3-e_4,v_6=e_1+e_2,v_7=e_1+e_2+e_3,v_8=-e_1-e_2$ & $V(v_1,v_4)$ & $-\displaystyle\frac{3}{2}$ \\ \hline

$I_{14}$ & $v_1=e_1, v_2=e_2,v_3=e_3,v_4=e_1+e_2,v_5=e_4,v_6=-e_1-e_2-e_3-e_4,v_7=-e_1-e_2-e_4,v_8=e_1+e_2+e_3+e_4$ & $V(v_1,v_4)$ & $-\displaystyle\frac{3}{2}$ \\ \hline

\end{tabular}

\begin{tabular}{  |c| p{11,2cm} | c| c| }
\hline
& \centering{Primitive \ Vectors} & Surface & $ch_2(T_X)\cdot S$  \\ 
\hline
$I_{15}$ & $v_1=e_1, v_2=e_2,v_3=e_3,v_4=e_4,v_5=-2e_1-2e_2-e_3-e_4,v_6=e_1+e_2,v_7=e_1+e_2+e_3,v_8=-e_1-e_2$ & $V(v_1,v_4)$ & $-\displaystyle\frac{3}{2}$\\ \hline
$J_1$ & $v_1=e_1, v_2=e_2,v_3=e_3,v_4=e_4,v_5=-e_3-e_4,v_6=e_1+e_2+e_3,v_7=e_1+e_2+2e_3,v_8=-e_1-e_2-e_3$ & $V(v_1,v_3)$ & -1 \\ \hline

$J_2$ & $v_1=e_1, v_2=e_2,v_3=-e_1-e_2,v_4=e_4,v_5=e_3,v_6=e_1+e_2-e_3-e_4,v_7=-e_3-e_4,v_8=-e_1-e_2+e_3+e_4$ & $V(v_1,v_3)$ & $-\displaystyle\frac{1}{2}$ \\ \hline
$K_1$ & $v_1=e_1, v_2=e_2,v_3=e_3,v_4=e_4,v_5=2e_1+2e_2-e_3-e_4,v_6=-e_1,v_7=-e_2,v_8=-e_1-e_2,v_9=e_1+e_2$ & $V(v_3,v_4)$ & -3 \\ \hline

$K_2$ & $v_1=e_1, v_2=e_2,v_3=e_3,v_4=e_4,v_5=2e_1+e_2-e_3-e_4,v_6=-e_1,v_7=-e_2,v_8=-e_1-e_2,v_9=e_1+e_2$ & $V(v_3,v_4)$ & -3 \\ \hline

$K_3$ & $v_1=e_1, v_2=e_2,v_3=e_3,v_4=e_4,v_5=e_1+e_2-e_3-e_4,v_6=-e_1,v_7=-e_2,v_8=-e_1-e_2,v_9=e_1+e_2$ & $V(v_3,v_4)$ & -3 \\ \hline

$M_1$ & $v_1=e_1, v_2=e_2,v_3=e_3,v_4=e_4,v_5=-e_4,v_6=e_1+e_2+e_3-e_4,v_7=-e_1-e_2-e_3+e_4,v_8=-e_1$ & $V(v_2,v_4)$ & $-\displaystyle\frac{5}{2}$ \\ \hline

$M_2$ & $v_1=e_1, v_2=e_2,v_3=e_3,v_4=e_4,v_5=e_1-e_4,v_6=e_1+e_2+e_3-e_4,v_7=-e_2-e_3+e_4,v_8=-e_1$ & $V(v_2,v_4)$ & $-\displaystyle\frac{5}{2}$ \\ \hline

$M_3$ & $v_1=e_1, v_2=e_2,v_3=e_3,v_4=e_4,v_5=e_1-e_4,v_6=e_1+e_2+e_3-e_4,v_7=-e_2-e_3,v_8=-e_1$ & $V(v_2,v_4)$ & $-\displaystyle\frac{5}{2}$ \\ \hline
$M_4$ & $v_1=e_1, v_2=e_2,v_3=e_3,v_4=e_4,v_5=e_1-e_4,v_6=e_1+e_2+e_3-e_4,v_7=-e_1-e_2-e_3+e_4,v_8=-e_1$ & $V(v_2,v_4)$ & $-\displaystyle\frac{5}{2}$\\ \hline

$M_5$ & $v_1=e_1,$ $v_2=e_2,$ $v_3=-e_1-e_2+e_4,$ $v_4=e_3,$ $v_5=e_1-e_3-e_4,$ $v_6=e_4,$ $v_7=e_1-e_4,$ $v_8=-e_3-e_4$ & $V(v_2,v_4)$ & $-\displaystyle\frac{3}{2}$ \\ \hline 

$Q_1$ & $v_1=e_1, v_2=e_2,v_3=e_3,v_4=e_4,v_5=e_1-e_3,v_6=e_1-e_4,v_7=-e_2,v_8=-e_1-e_2,v_9=e_1+e_2$ & $V(v_3,v_4)$ & $-\displaystyle\frac{3}{2}$ \\ \hline

$Q_2$ & $v_1=e_1, v_2=e_2,v_3=e_3,v_4=e_4,v_5=e_1-e_3,v_6=e_3-e_4,v_7=-e_2,v_8=-e_1-e_2,v_9=e_1+e_2$ & $V(v_3,v_4)$ & $-\displaystyle\frac{3}{2}$ \\ \hline

$Q_3$ & $v_1=e_1, v_2=e_2,v_3=e_3,v_4=e_4,v_5=e_1+e_2-e_3,v_6=e_1+e_2-e_4,v_7=-e_2,v_8=-e_1-e_2,v_9=e_1+e_2$ & $V(v_3,v_4)$ & $-\displaystyle\frac{3}{2}$ \\ \hline

$Q_4$ & $v_1=e_1, v_2=e_2,v_3=e_3,v_4=e_4,v_5=e_1-e_3,v_6=e_1+e_2-e_4,v_7=-e_2,v_8=-e_1-e_2,v_9=e_1+e_2$ & $V(v_3,v_4)$ & $-\displaystyle\frac{3}{2}$ \\ \hline

$Q_5$ & $v_1=e_1, v_2=e_2,v_3=e_3,v_4=e_4,v_5=e_1+e_2-e_3,v_6=e_3-e_4,v_7=-e_2,v_8=-e_1-e_2,v_9=e_1+e_2$ & $V(v_3,v_4)$ & $-\displaystyle\frac{3}{2}$ \\ \hline

$Q_7$ & $v_1=e_1, v_2=e_2,v_3=e_3,v_4=e_4,v_5=e_1+e_2-e_3,v_6=-e_2-e_4,v_7=-e_2,v_8=-e_1-e_2,v_9=e_1+e_2$ & $V(v_3,v_4)$ & $-\displaystyle\frac{3}{2}$ \\ \hline

$Q_9$ & $v_1=e_1, v_2=e_2,v_3=e_3,v_4=e_4,v_5=e_1+e_2-e_3,v_6=e_2-e_4,v_7=-e_2,v_8=-e_1-e_2,v_9=e_1+e_2$ & $V(v_3,v_4)$ & $-\displaystyle\frac{3}{2}$ \\ \hline

$Q_{12}$ & $v_1=e_1, v_2=e_2,v_3=e_3,v_4=e_4,v_5=e_1-e_3,v_6=e_2-e_4,v_7=-e_2,v_8=-e_1-e_2,v_9=e_1+e_2$ & $V(v_3,v_4)$ & $-\displaystyle\frac{3}{2}$ \\ \hline
\end{tabular}

\begin{tabular}{  |c| p{11,2cm} | c| c| }
\hline
& \centering{Primitive \ Vectors} & Surface & $ch_2(T_X)\cdot S$  \\ 
\hline
$Q_{13}$ & $v_1=e_1, v_2=e_2,v_3=e_3,v_4=e_4,v_5=e_2-e_3,v_6=e_2-e_4,v_7=-e_2,v_8=-e_1-e_2,v_9=e_1+e_2$ & $V(v_3,v_4)$ & $-\displaystyle\frac{3}{2}$ \\ \hline

$Q_{14}$ & $v_1=e_1, v_2=e_2,v_3=e_3,v_4=e_4,v_5=e_2-e_3,v_6=e_3-e_4,v_7=-e_2,v_8=-e_1-e_2,v_9=e_1+e_2$ & $V(v_3,v_4)$ & $-\displaystyle\frac{3}{2}$ \\ \hline

$Q_{16}$ & $v_1=e_1, v_2=e_2,v_3=e_3,v_4=e_4,v_5=e_1+e_2-e_3,v_6=-e_1-e_2-e_4,v_7=-e_2,v_8=-e_1-e_2,v_9=e_1+e_2$ & $V(v_3,v_4)$ & $-\displaystyle\frac{3}{2}$ \\ 
\hline

$Q_{17}$ & $v_1=e_1, v_2=e_2,v_3=e_3,v_4=e_4,v_5=e_2-e_3,v_6=-e_1-e_2-e_4,v_7=-e_2,v_8=-e_1-e_2,v_9=e_1+e_2$ & $V(v_3,v_4)$ & $-\displaystyle\frac{3}{2}$ \\ \hline

$R_1$ & $v_1=e_1, v_2=e_2,v_3=e_3,v_4=-e_1-e_2+e_3,v_5=-e_1-e_2,v_6=e_4,v_7=-e_4,v_8=e_1+e_2-e_3-e_4,v_9=-e_1-e_2+e_3+e_4$ & $V(v_1,v_3)$ & -4 \\ \hline

$R_2$ & $v_1=e_1, v_2=e_2,v_3=e_3,v_4=-e_1-e_2+e_3,v_5=-e_1-e_2+e_4,v_6=e_4,v_7=-e_4,v_8=e_1+e_2-e_3-e_4,v_9=-e_1-e_2+e_3+e_4$ & $V(v_1,v_3)$ & -4 \\ \hline
$R_3$ & $v_1=e_1, v_2=e_2,v_3=e_3,v_4=-e_1-e_2+e_3, v_5=-e_3, v_6=e_4, v_7=-e_4,v_8=e_1+e_2-e_3-e_4,v_9=-e_1-e_2+e_3+e_4$ & $V(v_1,v_3)$ & -4 \\ \hline

$108$ & $v_1=e_1, v_2=e_2,v_3=e_3, v_4=-e_1-e_2+e_4,v_5=-e_1-e_2-e_3+e_4,v_6=-e_3,v_7=-e_4,v_8=e_1-e_4,v_9=e_4$ & $V(v_4,v_9)$ & -1 \\ \hline

$U_1$ & $v_1=e_1, v_2=e_1+e_3,v_3=e_3,v_4=-e_1,v_5=-e_1-e_3,v_6=-e_3,v_7=e_2,v_8=e_1-e_2,v_9=e_4,v_{10}=e_1-e_4$ & $V(v_3,v_7)$ & $-\displaystyle\frac{1}{2}$ \\ \hline

$U_2$ & $v_1=e_1, v_2=e_1+e_3,v_3=e_3,v_4=-e_1,v_5=-e_1-e_3,v_6=-e_3,v_7=e_2,v_8=e_1-e_2,v_9=e_4,v_{10}=e_1-e_2-e_4$ & $V(v_3,v_7)$ & $-\displaystyle\frac{1}{2}$ \\ \hline
$U_3$ & $v_1=e_1, v_2=e_1+e_3,v_3=e_3,v_4=-e_1,v_5=-e_1-e_3,v_6=-e_3,v_7=e_2,v_8=e_1-e_2,v_9=e_4,v_{10}=e_1+e_3-e_4$ & $V(v_3,v_9)$ & $-\displaystyle\frac{1}{2}$ \\ \hline
$U_7$ & $v_1=e_1, v_2=e_1+e_3,v_3=e_3,v_4=-e_1,v_5=-e_1-e_3,v_6=-e_3,v_7=e_2,v_8=e_1-e_2,v_9=e_4,v_{10}=e_3-e_4$ & $V(v_3,v_9)$ & $-\displaystyle\frac{1}{2}$ \\ \hline
$U_8$ & $v_1=e_1, v_2=e_1+e_3,v_3=e_3,v_4=-e_1,v_5=-e_1-e_3,v_6=-e_3,v_7=e_2,v_8=e_1-e_2,v_9=e_4,v_{10}=-e_1-e_4$ & $V(v_3,v_9)$ & $-\displaystyle\frac{1}{2}$ \\ \hline
$Z_1$ & $v_1=e_1, v_2=e_2,v_3=e_3,v_4=e_4,v_5=-e_1-e_2,v_6=-e_2-e_3-e_4,v_7=e_1-e_3-e_4,v_8=-e_1+e_4$ & $V(v_1,v_3)$ & $-\displaystyle\frac{5}{2}$ \\ \hline

$Z_2$ & $v_1=e_1, v_2=e_2,v_3=e_3,v_4=e_4,v_5=-e_1-e_2,v_6=-e_3-e_4,v_7=e_1+e_2-e_3-e_4,v_8=-e_1+e_4$ & $V(v_1,v_3)$ & -2 \\ \hline

$117$ & $v_1=e_2, v_2=e_3,v_3=e_4,v_4=-e_1,v_5=-e_2,v_6=-e_3,v_7=-e_4,v_8=e_1+e_2+e_3+e_4,v_9=-e_1-e_2-e_3-e_4,v_{10}=e_1$ & $V(v_1,v_4)$ & -5 \\ \hline

$118$ & $v_1=e_2, v_2=e_3,v_3=e_4,v_4=-e_1,v_5=-e_2,v_6=-e_3,v_7=-e_4,v_8=e_1+e_2+e_3+e_4,v_9=e_1,$ & $V(v_1,v_4)$ & $-\displaystyle\frac{5}{2}$ \\ \hline

$124$ & $v_1=e_1, v_2=e_2, v_3=-e_1-e_2,v_4=-e_1+e_4, v_5=e_1-e_3-e_4, v_6=e_3, v_7=e_4, v_8=e_1+e_2-e_3-e_4, v_9=-e_1-e_2+e_3,$ & $V(v_1,v_7)$ & $-4$ \\ \hline
\end{tabular}\\ \\

\textbf{Remark}. We have concluded that there exists only one toric 2-Fano 4-fold. However, using proposition 3.2 and computations with the program Maple we see that there exist toric Fano 4-folds that have nef second Chern character (i.e., $ch_2(T_X)\cdot S\geq 0$ for every surface $S\subset X$). They are: \\$\mathbb{P}^4, B_1, B_2, B_3, B_4, C_4, D_1, D_2, D_3, D_5, D_6, D_8, D_9, D_{12}, D_{13}, D_{15}, L_1, L_2, L_3, L_4, L_5, L_6, L_7, L_8, L_9$.

\section{Appendix: Maple Code}

In this appendix we provide the code used in the program Maple to compute $ch_2(T_X)\cdot S$ for toric 4-folds.
The reader who wishes to obtain the file in Maple extension can access the webpage \url{http://w3.impa.br/~edilaine/}

\begin{verbatim}
restart:
with(LinearAlgebra):

v[0]:=<0,0,0,0>:
# Insert the minimal vectors of a toric 4-fold.
v[1]:=< >:
# Insert each primitive collection that consists of two vectors {v_i,v_j} by {i,j}.
A[1]:={}:
# Insert each primitive collection that consists of three vectors {v_i,v_j,v_k} by {i,j,k}.
B[1]:={}:
# Insert each primitive collection that consists of four vectors {v_i,v_j,v_k,v_r} by {i,j,k,r}.
C[1]:={}:
# Insert the number of minimal vectors.
z:=
# Insert the number of primitive collections that consist of two vectors.
a:=
# Insert the number of primitive collections that consist of three vectors.
b:=
# Insert the number of primitive collections that consist of four vectors.
c:=
j:=1:
f:=(x,y,z)->if(x=y or x=z,0,x):
g:=(x,y,z)->if(x=y,z,y):
h:=(x,y,z,k)->if(x=y,<z,k>,if(x=z,<y,k>,<y,z>)):
while j<z+1 do
for i from 1 to z do
if {i,j} in [seq(A[k],k=1..a)]or i=j then F[i,j]:=[0,0,0] else F[i,j]:=[i,j,0] end if:
Q[0][i,j]:=[0,0,0]: end do: j:=j+1 end do:
for w from 1 to z do
for i from 1 to z do
for j from 1 to z do
V[i,j,w]:=if({i,j,w} in [seq(B[k], k=1..b)] or i=j or j=w or i=w or F[i,j]=[0,0,0] 
or F[i,w]=[0,0,0] or F[j,w]=[0,0,0], [0,0,0], [i,j,w]):
if (w!=i and w!=j or F[i,j]=[0,0,0]) then Q[w][i,j]:=V[i,j,w] else 
K[w,i,j]:=v[w]-Norm(v[g(w,i,j)],2)^(-2) * DotProduct(v[w],v[g(w,i,j)])* v[g(w,i,j)]:
Q[w][i,j]:=seq([Q[f(k,i,j)][i,j],-(DotProduct(v[w],K[w,i,j]))^(-1) * DotProduct(v[f(k,i,j)],
K[w,i,j])],k=1..z) end if end do end do end do
H:=(i,j,m,w)->if({i,j,m,w} in [seq(C[k],k=1..c)] or V[i,j,m]=[0,0,0] or V[i,j,w]=[0,0,0]
or V[i,m,w]=[0,0,0] or V[w,j,m]=[0,0,0],0,1):
for w from 1 to z do
for i from 1 to z do
for j from 1 to z do
for m from 1 to z do
P[w][0,0,0]:=0:
P[0][i,j,m]:=0:
if (w!=i and w!=j and w!=m) then P[w][i,j,m]:=H(i,j,m,w) else if V[i,j,m]=[0,0,0] then
P[w][i,j,m]:=0 else q[w,i,j,m]:=GramSchmidt([v[h(w,i,j,m)[1]],v[h(w,i,j,m)[2]],v[w]):
u[w,i,j,m]:=-q[w,i,j,m][3]*DotProduct(v[w],q[w,i,j,m][3])^(-1):
P[w][i,j,m]:=add(DotProduct(u[w,i,j,m],v[n])*H(i,j,m,n),n=1..z):
end if:
end if end do end do end do end do
for w from 1 to z do
for i from 1 to z do
for j from 1 to z do
for r from 1 to z do
C[r,w,i,j]:=if(Q[w][i,j]=[0,0,0],0,if(Q[w][i,j]=[i,j,w],P[w][i,j,w],P[w][Q[w][i,j][r][1][1],
Q[w][i,j][r][1][2],Q[w][i,j][r][1][3]]*Q[w][i,j][r][2])):
M[w][i,j]:=if(Q[w][i,j]=[0,0,0],0,if(Q[w][i,j]=[i,j,w],P[w][i,j,w],add(C[p,w,i,j],p=1..z))):
end do end do end do end do
Y:=(a,b)->seq(M[j][a,b],j=1..z):
for x from 1 to z do
for y from x+1 to z do
ch2[x,y]:=add(M[j][x,y],j=1..z)/2:
#ch2[i,j] is the intersection number ch_2(X).V(v_i,v_j). 
if evalb({x,y} in [seq(A[k],k=1..a)])=false then print([x,y],ch2[x,y]) else end if;
end do end do:
\end{verbatim}

\textbf{Acknowledgments}\\

I wish to thank my advisor, Carolina Araujo, who pointed out to me Batyrev's paper on Fano toric 4-folds \cite{bat}, for suggesting this question, and for constantly guiding me. I am also grateful to Alex Abreu and Douglas Monsôres for interesting and useful discussions.

\end{document}